\documentclass[secthm,seceqn,amsthm,ussrhead,12pt]{amsart}
\usepackage[russian]{babel}
\usepackage{amsmath,latexsym}
\usepackage[psamsfonts]{amssymb}
\usepackage{graphicx}
\usepackage{caption2}
\usepackage{times}
\usepackage[mathcal]{euscript}
\numberwithin{equation}{section} \textwidth=16cm \topmargin=0.1cm
\renewcommand{\epsilon}{\varepsilon}

\newcommand{\be}{\begin{equation}}
\newcommand{\ee}{\end{equation}}

{\bf}{\it}

\begin{document}
УДК 517.518.5

\author{ Д.Д.Туракулов $^1$}
\address{$^1$ Туракулов Давир Давронович.   г Самарканд:
Раб. адрес., 703004, Университетский бульвар-15, СамГУ
механика-математический факультет, кафедра "Мате\-ма\-ти\-чес\-кая
физика и теория функций". E-mail: davirt@rambler.ru}

\title[Об ограничении преобразования Фурье на
гиперповерхностях] {Об ограничении преобразования Фурье на
гиперповерхностях} \maketitle

{\bf Аннотация.} Рассматривается преобразования  Фурье на выпуклых
аналитических гиперповерхностях в $R^{4}.$

Получено решение проблемы об $\left(L_{p}, L_{2} \right)$
ограниченности соответствующего оператора ограничения
преобразования Фурье на некоторых классов гиперповерхностей.

{\bf Ключевые слова:} преобразование Фурье, гиперповерхность
конечного типа, оператор ограничения, показатель осцилляции,
высота функции.

\begin{center}
\section{Введение}
\end{center}

Как известно, преобразование Фурье суммируемой в $R^{n} $ функции
является  непре-рывной функцией. Поэтому мы можем определить ее в
каждой точке двойственного прост-ранства. С другой стороны теорема
Планшареля показывает, что  преобразо-вание Фурье $L_{2}
(R_{x}^{n})$ функции снова является $L_{2} (R_{\xi }^{n})$
функцией и обратно. Следова-тельно, преобразование Фурье $L_{2}
(R_{x}^{n})$ функции может быть любой функцией из $L_{2}(R_{\xi
}^{n}).$ Более того неравенство Хаусдорфа-Юнга показывает, что
преобразование Фурье $L_{p} (R_{x}^{n})$ функции при  $1\le p\le
2$ может быть определено как функция из класса $L_{q} (R_{\xi
}^{n}),$ где $q$- двойственное число, т.е. $\frac{1}{p}
+\frac{1}{q}=1.$ Несмотря на то, что преобразование Фурье в $R^{n}
$ определяется линейными структурами $R^{n},$ результаты об
ограничении преобразования Фурье пока-зывают связь между
гармоническим анализом и нелиней-ным объектом. Точнее, как показал
Фефферман [8] преобразование Фурье $L_{p}$ функции при достаточно
близких к единице корректно определена как $L_{2} \left(S^{n-1}
\right)$ (где $S^{n-1}$- единичная сфера в $R^{n}$) функции, хотя
$S^{n-1} $ имеет меру нуль в $R^{n} $ и $\hat{f}$ априори
определено почти всюду. Здесь можно отметить классические
результаты Томаса [5], Штрихартса [3], Стейна [4], Феффермана [8],
Зигмунда [6].

Аналогичные проблемы возникают в случае, когда $S$ произвольная
поверхность Евклидова пространства.

Перейдем к более точным определениям.

Пусть $f\in Sh\left(R^{n} \right),$ где $Sh\left(R^{n}\right)$
-пространство Шварца. Через $S\subset R^{n} $ обозначим гладкую
поверхность в пространстве $R^{n} $ и $dS$-элемента индуцированной
лебеговой меры, а также введем меру $d\mu (x)=\psi (x)dS,$ где
$\psi \in C_{0}^{\infty } (S)$ неотрицательная гладкая функция с
компактным носителем. Пространство $L_{q} (S)$ определяется
естественным образом.

Если $f\in Sh\left(R^{n} \right),$ то естественно  определен
оператор   $(Rf)(\xi)=\hat{f}(\xi),$ при $\xi \in S,$ так как
$\hat{f}(\xi)$ является элементом пространства Шварца
$Sh\left(R_{\xi}^{n} \right)$.

\textbf{Определение 1.1.} {\it Говорят, что оператор $R$ имеет тип
$\left(L_{p}, L_{q} \right)$ если существует положительное число
$A_{p,q} $ такое, что для любой функции  $f\in Sh \left(R^{n}
\right)$  выполняется неравенство:
\begin{equation} \label{1}
\left(\int _{S}\left|\hat{f}\right|^{q} d\mu  \right)^{\frac{1}{q}
} \le A_{p,q} \left\| f\right\|_{L_{p} }
\end{equation}
\textit{где $\left\| f\right\| _{L_{p} } =\left(\int_{R^{n}
}\left|f(x)\right|^{p} dx \right)^{\frac{1}{p}}.$}

Проблема о точных значениях $p,q$ до сих пор остается открытой.
Эта проблема называется "проблемой об ограничении преобразования
Фурье".

Классические  результаты Томаса и Стейна [4] показывают,  что если
$S$ единичная сфера и  $q=2$, то при $ 1\le p\le
\frac{2(n+1)}{n+3} $ имеет  место неравенство вида \eqref{1} для
любой функции   $f\in Sh\left(R^{n} \right).$

 Мы рассмотрим эту задачу в случае, когда $q=2$ и $S$ выпуклая
 аналитическая гиперповерхность в $R^{4}.$
 Локально мы можем определить $S$ как график аналитической функции

\[x_{4} =\Phi(x_{1} ,x_{2} ,x_{3} ),\]
 удовлетворяющая условиям:
 $\Phi(0)=0,{\rm \; \; \; \; \; \; \; }\nabla \Phi (0)=0$.

В работе [10] рассмотрена аналогичная задача и получены  некоторые
оценки для случая, когда $S$ удовлетворяет некоторым условиям, так
называемых условиям полиэд-рального типа.

 Нам необходимо некоторые обозначения и определения для
 формулировки основного результата.

Пусть $f$ гладкая функция в некоторой окрестности нуля и она
удовлетворяет  условиям:

\[f(0)=0,{\rm \; \; \; }\nabla f(0)=0.\]
Функция $f$ называется \textit{выпуклой,} если для любых векторов
$x,y\in U\subset R^{n} $  (где $U$-некоторая выпуклая окрестность
начала координат) и для любых неотрицательных чисел $\alpha ,\beta
,$ удовлетворяющих условию $\alpha +\beta =1$ имеет место
неравенство:

\[f(\alpha x+\beta y)\le \alpha f(x)+\beta f(y).\]
Функция $f$ называется \textit{функцией конечного линейного типа}
в  начале координат, если для любого единичного вектора $\xi \in
R^{n} $ существует $N\ge 2$ такое, что $D_{\xi }^{N} f(0)\ne 0$,
где $D_{\xi }^{N} f(0)$ производная функции $f$ по направлению
вектора $\xi $ в начале координат. Поверхность $S$ называется
\textit{выпуклой конечного линейного типа, }если она локально
задается в виде графика выпуклой функции конечного линейного типа.

Пусть $f$ гладкая функция, определенная в начале координат.
Рассмотрим  ее ряд Тейлора с центром в начале координат(т.е. ряд
Маклорена этой функции)

\[f_{x} \approx \sum _{k\in N_{0}^{n} }c_{k} x^{k}  ,{\rm \; \; \;
\; \;  \; \; \; }c_{k} \in R.\] Носитель (носителем Тейлора) этого
ряда определяется соотношением:

 $$
\tau {\rm \; }(f_{x})=\left\{k\in N_{0}^{n} \backslash \{ 0\}
:c_{k} \ne 0\right\}.$$

 Многогранник Ньютона ряда Маклорена $f$ определен как выпуклая оболочка
 мно-жества  $\bigcup \left\{k+R_{+}^{n} \right\}$, где $k\in \tau {\rm \;
  }(f_{x} )$  и $R_{+}^{n} =\left(R_{+} \right)^{n} $.

Фиксируем систему координат в $R^{n} $ и обозначим через $f_{x} $
ряд  Маклорена функции $f$ в этой системе координат. Пусть $d$
координата пересечения прямой $x_{1} =\ldots =x_{n} =t,t\in R$ с
границей многогранника Ньютона. Это число будет называться
\textit{расстоянием }между многогранником и началом координат.
Расстояние обозначается через $d(x)$. \textit{Главной гранью}
многогранника Ньютона называется грань минимальной размерности,
содержащей точку $(d(x),\ldots ,d(x))$.

Пусть  $f$ функция, определенная выше, и пусть $x=(x_{1} ,\ldots
,x_{n} ),$ фиксированная система координат в нуле в $R^{n} $.
Обозначим через $f_{x} $ ряд Маклорена $f$, $d(x)$ расстояние
между началом координат и многогранником Ньютона $N(f_{x} )$.
Рассмотрим величину $h(f)=\sup \left\{d(x)\right\}$, где
"supremum" берется относительно набора всех локальных глад-ких
систем координат $x$ в начале координат. Число $h(f)$ называется
\textit{высотой  }функции $f$. Варченко А.Н доказал существование
так называемых "приспособленных" систем координат, т.е таких
систем координат, где выполняется равенство: $h(f)=d$.
Существо-вание аналогичных систем координат для гладких выпуклых
функций конечного линейного типа доказано в работе [7] и для
произвольных выпуклых аналитических функций в работе [11].

Рассмотрим меру $d\mu (x)=\psi (x)dS$ и преобразование Фурье этой
меры :

\[d\hat{\mu }(\xi )=J(\xi ):={\rm \; }\int _{S}\exp (i(x,\xi ))d\mu (x) .\]
В работе [12] доказана следующая

\textbf{Теорема 1.1.} \textit{Пусть $\Phi $ выпуклая аналитическая
функция в окрестности нуля. Тогда существует окрестность нуля  $U$
такая, что для любой функции  $\psi \in C_{0}^{\infty } (U)$
справедлива следующая оценка}

\[\left|\mathop{d\mu (\xi )}\limits^{\wedge }
 \right|\le \frac{C\left\| \psi \right\| _{C^{(4)} (U)} }{\left|\xi
 \right|^{\frac{1}{h(\Phi )} } } ,\]
 где  $h(\Phi )$  высота функции $\Phi $ , $\left\| \psi
\right\| _{C^{(4)} (U)} $ норма функции на $C^{(4)} $.

Основным результатом настоящей работы является следующая

\textbf{Теорема 1.2.   }\textit{Пусть  }$\Phi$ \textit{-выпуклая
аналитическая функция в окрестности нуля в $R^{3} $. Тогда
существует окрестность нуля $U$ такая, что для любой функции
$\psi \in C_{0}^{\infty } (U)$ оператор ограничения $R$ имеет тип
$\left(L_{p} ,L_{2} \right)$ при любых $1\le p\le \frac{2(1+h(\Phi
))}{2h(\Phi )+1} $. Более того, если  $\psi (0)>0$, то оператор
ограничения имеет тип $\left(L_{p} ,L_{2} \right)$ тогда и только
тогда, когда  $1\le p\le \frac{2(1+h(\Phi ))}{2h(\Phi )+1} $.}

\section {Оценка сверху для оператора ограничения.}

\textbf{ Доказательство} результата предыдущего пункта
основывается на следущей теоремы А.Гринлифа [2] и теореме 1.1 о
равномерной оценке соответствующего осцилляторного интеграла.

\textbf{Теорема 2.1.} \textit{Пусть $S\subset R^{n+m} $ гладкое
$n$-мерное подмногообразие с гладкой мерой $d\mu =\psi (x)dS$,
сосредоточенной вне границы $S$. Предположим для некоторых $C,r>0$
выполняется неравенство}

\[\left|d\hat{\mu }(\xi )\right|\le C\left(1+\left|\xi \right|\right)^{-\beta } .\]
\textit{Тогда существует $C'>0$ такое, что для любой  функции
$f\in Sh\left(R^{n+m} \right)$ выполняется неравенство}

\begin{equation} \label{GrindEQ__2_}
\left\| \left. \hat{f}\right|_{S} \right\| _{L_{2} (S,d\mu )}
\le C'\left\| f\right\| _{L_{p} (R^{n+m} ,dx)} ,
\end{equation}
\textit{при $p=\frac{2(m+\beta )}{2m+\beta } $.}

В частности если $S$ гладкая гиперповерхность,  то $m=1$ и мы
имеем оценку с $ $\textit{$p=\frac{2(1+\beta )}{2+\beta } $. }Рады
полноты изложения мы остановимся в основных моментах
доказательства теоремы 2.1 А.Гринлифа при  $m=1$ [2].

Сначала заметим, что оценка \eqref{GrindEQ__2_} инвариантно
относительно Евклидова движения пространств $R_{x}^{4} $ и $R_{\xi
}^{4} $, так как $dS$ инвариантно относительно таких
преобразований. Следовательно, мы можем считать что $S$ задана в
виде графика гладкой функции $x_{4} =\Phi (x_{1} ,x_{2} ,x_{3} )$,
удовлетворяющей условиям $\Phi (0)=0,{\rm \; \; \; \; \; }\nabla
\Phi (0)=0$. При этом используя соответствующего разбиения единицы
мы можем считать, что носитель функции $\psi $ сосредоточен в
достаточно малой  окрестности начало координат.

Теперь с помощью леммы Томаса [5]  достаточно показать
справедливость оценки

\[\left\| d\hat{\mu }*g\right\| _{L_{q} } \le C\left\| g\right\| _{L_{p} } ,\]
где  $*$- свертка, $p=\frac{2(1+\beta )}{2+\beta } $ и $\frac{1}{p} +\frac{1}{q} =1$.

Рассмотрим обобщенную функцию (семейство аналитических обобщенных функций)

\[G_{z} (x)=\psi (x_{1} ,x_{2} ,x_{3} ,x_{4} )
\frac{\left|x_{4} -\Phi (x_{1} ,x_{2} ,x_{3} \right|^{z}
}{\Gamma(z+1)} .\] Как известно, $G_{z} $ аналитически
продолжается в $C$ [1] и

\[G_{-1} =\psi \delta (x_{4} -\Phi (x)),\]
где $\delta (x_{4} -\Phi (x))$- "дельта" функция на поверхности
$S$.

Теперь определим семейство аналитических операторов

\[T_{z} f(x)=\hat{G}_{z} *f(x),\]
где $f\in Sh\left(R^{n} \right)$.

Покажем справедливость неравенство: $\left|\frac{1}{\Gamma(1+iy)}
\right|\le Ce^{\pi \left|y\right|} $.

Заметим, что $\Gamma(1 )=1$. Поэтому  из аналитичности $\Gamma $ в
точке $1$ имеем $\left|\Gamma (z)\right|\ge \frac{1}{2} $ при
$\left|z-1\right|<\delta $ для некоторого положительного число
$\delta $. В частности при $\left|y\right|<\delta$  мы имеем
$\frac{1}{\left|\Gamma(1+iy)\right|} \le 2$.

Теперь используем соотношение  $\Gamma(z)\Gamma(1-z)=\frac{\pi
}{\sin \pi z} $.

При $z=-iy$ мы получим
\[\frac{\pi }{\Gamma(1+iy)} =\sin \pi (1+iy)\Gamma(-iy)=
-\sin \pi iy\Gamma(-iy)=ish\pi y\Gamma(-iy).\] Отметим, что если
$\left|y\right|\ge \delta $, то существует $C_{\delta } $  такое,
что  $\left|\Gamma(-iy)\right|\le C_{\delta } $.

Действительно, по формуле аналитического  продолжения гамма
функции, имеем следующие соотношения:
\[\left|\Gamma(1-iy)\right|\le \frac{\left|\Gamma(1-iy)\right|}{\left|iy\right|}
 \le \frac{1}{\delta } .\]
Итак, $\frac{1}{\left|\Gamma(1+iy)\right|} \le C_{\delta } e^{\pi
\left|y\right|} $.

Если $Rez=0$, то  $\mathop{\sup }\limits_{x} \left|G_{z}
(x)\right|=C_{z} \le Ce^{\pi \left|Imz\right|} $.$ $

Следовательно,

\[\left\| T_{z} g\right\| _{L_{2} } \le C_{z} \left\| g\right\| _{L_{2} } .\]
Теперь используем следующую формулу [1]:
\[\hat{x}^{z} =2^{z+1} \sqrt{\pi } \frac{\Gamma(z+1)}{\Gamma(-z)}
\left|\xi \right|^{-z-1} .\] Поэтому
\[\hat{G}_{z} (\xi )=\frac{2^{z+1} \sqrt{\pi } \left|\xi _{4}
\right|^{-z-1} }{\Gamma(-z)} d\hat{\mu }(\xi ).\] Отсюда вытекает,
что для  $Rez=-\beta -1$ выполняется следующее неравенство

$\left|\hat{G}_{z} (\xi )\right|\le C_{z} $, причем $\left|C_{z}
\right|\le Ce^{C_{1} \left|Imz\right|} $. Следовательно,

\[\left\| T_{z} g\right\| _{L^{\infty } } \le C_{z} \left\| g\right\| _{L_{1} } .\]
Используя интерполяционную теорему при $p=\frac{2(1+\beta
)}{2+\beta } $ мы имеем:

\[\left\| T_{-1} g\right\| _{L_{q} } \le C_{k} \left\| g\right\| _{L_{p} } \]
Наконец, используя теоремы 1.1 мы придем к оценку сверху для
оператора ограничения.

\section{Оценка снизу для оператора ограничения.}

\textbf{ }Как отметили выше $\left(L_{p} ,L_{2} \right)$
ограниченность оператора сужения $R$ инвариантно относительно
аффинного преобразования. В частности мы можем считать, что
функ-ция $\Phi $ записана в координате Шульца [7], в случае когда
$S$ выпуклая гиперповерхность.

Обозначим через $d$ расстояние до многогранника Ньютона функции
$\Phi $  построенного в точке $x=0$. Следующее утверждение
доказано в работе [9]:

\textbf{Лемма 3.1. } \textit{Пусть  $q$ сопряженное число, т.е
$\frac{1}{q} +\frac{1}{p} =1{\rm \; \; }(p\ge 1)$. Если $\psi
(0)>0$ и оператор ограничения имеет тип }$\left(L_{p} ,L_{2}
\right)$, \textit{то для сопряженного числа $q$ справедливо
неравенство $q\ge 2d+2$. Если $S$ выпуклая поверхность заданная в
виде графика функции $\Phi $ и  $N(\Phi )$  многогранник Ньютона
функции $\Phi$, то существует линейное преобразование такое, что
расстояние } $d$ \textit{определяется из соотношения $d=h(\Phi)$.}

\textbf{Доказательство теоремы 1.2.} Так как свойство
ограниченности оператора $R$ инвари-антно относительно линейных
преобразований, то из леммы 3.1 следует, что $1\le p\le
\frac{2(1+h(\Phi ))}{2h(\Phi )+1} $.  Если $q=2\left(1+h(\Phi
)\right)$, то
\[\frac{1}{p} =1-\frac{1}{q} =1-\frac{1}{2(h(\Phi )+1)} =
\frac{2h(\Phi )+1}{2(1+h(\Phi ))}.\] Отсюда вытекает
доказательство теоремы 1.2.

Автор приносит глубокую благодарность д.ф-м.н. И.А.Икромову за
постановки задачи и внимание к работе.

Работа финансирована грантом ОТ-Ф1-006 Госкомитета Науки и техники
Республики Узбекистан.

\vspace{1cm}

 \textbf{СПИСОК ЛИТЕРАТУРЫ}\textit{}

[1] Гельфанд И. М., Шилов Г. Е. Обобщенные функции и действия над
ними.  Вып. 1. М.: Физматгиз. 1959.

[2]  Greenleaf A. Principal curvature and Harmonic Analysis,
Indiana. Math. J. 30, 4, (1981), 519--537.

[3]   Srichartz R. Restrictions of Fourier transforms to quadratic
surfaces, Duke Math. J. 44, (1977), 705--714.

[4]  Stein E. M. Harmonic analysis: Real-valued methods,
Orthgonality, and oscillatory integrals. Princeton Univ. press,
43, 1993.

[5]  Tomas P. A. A restriction theorem for the Fourier transform,
Bull. Amer. Math. Soc. 81, (1975), 477--478.

[6]  Zygmund A. On Fourier coefficients and transforms of
functions of two variables. Studia Math. 50, (1974), 189--201.

[7]  Shulz H. Convex hupersufaces of  finite type and the
asymptotics of their Fourier transforms, Indiana Univ. Math. J.
40:4, (1991), 1267--1275.

[8]  Fefferman C. Inequalities for strongly singular convolution
operators, Acta  Math. 124, (1970), 9--36.

[9]  Magyar A. On Fourier restriction and Newton polygon.//Proc.
AMS, 137:2, (2009), 615--625.

[10] Iosevich A. Fourier transform, $L_{2} $ restriction theorem
and scaling// Boll.  Unione. Math. Italy, 8:2, (1999), 383--387.

[11] Икромов И.А и Солеев А.С. Алгоритм определения базиса Щульца
для выпуклых функций// Узбекский Мат. Журнал. 4, (2008), 75--88.

[12] Туракулов Д.Д.  Равномерные оценки осцилляторных интегралов с
выпуклой фазой в трехмерном случае// Докл. АН РУз. Ташкент, 1,
(2010), 7--12.

 \section*{On restriction of the Fourier transform to hypersufaces.}

 \centerline{D.D.Turakulov}

\textbf{Abstract.} It is considered Fourier transform of convex
analytic hypersufaces on $R^{4} $. We  prove that the
Fourier\textbf{ }restriction operator associated to convex
analytic hypersufaces is \textit{$\left(L_{p} ,L_{2} \right)$}
bounded whenever  $1\le p\le \frac{2h+2}{h+2} $. The result is
sharp.

\vspace{2cm}

\end{document}